\documentstyle[12pt]{article}
\input mssymb
\newtheorem{tw}{Theorem}
\newtheorem{df}{Definition}
\def\ot{\o \theta}
\def\ota{\ot_\ka}
\def\otb{\ot_\kb}
\def\d{\delta}
\def\g{\gamma}
\def\kd{{\dot\delta}}
\def\kg{{\dot\gamma}}
\def\graded{$Z_2$-graded }
\def\H{{\cal H}}

\def\1{{(1)}}
\def\2{{(2)}}
\def\[{\left[}
\def\]{\right]}

\def\i{\item}
\def\k{\kappa}
\def\s{\sigma}
\def\U{{\cal U}}

\def\trl{\triangleleft}
\def\ftt{\footnotetext}
\def\ftm{\footnotemark}
\def\bib{\bibitem}
\def\sec{\setcounter{equation}{0}}
\title{ Bicrossproduct Hopf Superalgebras and $D=4$ \kdef
\poin Supergroup}

\author{\em P.\ Kosi{\'n}ski \ftm[1]  \ftm[5]  ,
            J.\ Lukierski \ftm[2]  \ftm[6]  ,
            P.\ Ma{\'s}lanka \ftm[3]  \ftm[5],
             and J.\ Sobczyk \ftm[2]  \ftm[6]
}
\date{January 1995}

\def\<{\left<}
\def\>{\right>}

\def\e{\epsilon}
\def\ben{\begin{enumerate}}
\def\een{\end{enumerate}}
\def\:{\,\,:\,\,\,}

\def\to#1{\,\,{\stackrel{#1}\longrightarrow}\,\,}
\def\Pc{{\cal P}_4}
\def\({\left(}
\def\){\right)}
\def\tens{\otimes}
\def\~{\widetilde}

\def\poin{Poincar\'e }

\def\bowti{\triangleright\!\!\!<}

\def\ptc{>\!\!\!\blacktriangleleft}
\def\btc{\,\cobicross\,}
\def\ctb{\,\bicross\,}
\newcounter{popnr}
\renewcommand{\theequation}{\arabic{section}.\arabic{equation}}
\def\alpheqn{\setcounter{popnr}{\value{equation}}
             \stepcounter{popnr}
             \setcounter{equation}{0}
             \def\theequation{\arabic{section}.\arabic{popnr}\alph{equation}}
             }
\def\reseteqn{\setcounter{equation}{\value{popnr}}
              \def\theequation{\arabic{section}.\arabic{equation}}
              }
\newcommand{\beq}{\begin{eqnarray}}
\newcommand{\eeq}{\end{eqnarray}}
\newcommand{\beqq}{\begin{eqnarray*}}
\newcommand{\eeqq}{\end{eqnarray*}}

\def\bel#1{\begin{equation}\label{#1}}
\def\be{\begin{equation}}
\def\ee{\end{equation}}
\def\o{\overline}

\def\sP{{\cal P}_{4;1}}
\def\sT{T_{4;2}}

\def\r#1{(\ref{#1})}
\def\t{\tau}

\def\bl{\alpheqn}
\def\el{\reseteqn}
\def\a{\alpha}
\def\b{\beta}
\def\ka{{\dot\a}}
\def\kb{{\dot\b}}
\def\ba{\begin{array}}
\def\ea{\end{array}}
\def\cop{\Delta}

\def\so{SO(3,1;2)}
\def\kdef{$\k$-deformed }
\def\0{^{(0)}}

\begin{document}
\begin{titlepage}
%\rightline{ IFTUWr 887/95}
\maketitle
\thispagestyle{empty}
\begin{abstract}
The general framework of bicrossproduct Hopf algebras given by Majid is
extended to $Z_2$-graded bicrossproduct Hopf superalgebras.
As examples of bicrossproduct Hopf superalgebras we provide the graded
algebras of functions on undeformed as well as \kdef $D=4$ \poin
supergroups.
\end{abstract}
\vfill
\rightline{ IFTUWr 887/95, q-alg/9501010}
\def\thefootnote{\fnsymbol{footnote}}

\ftt[1]{Institute of Physics, University of
{\L}\'od\'z, ul. Pomorska 149/153, 90-236 {\L}\'od\'z, Poland.}
\ftt[2]{
Institute for Theoretical Physics, University of Wroc{\l}aw,
pl. Maxa Borna 9, 50-204 Wroc{\l}aw, Poland.}
\ftt[3]{Dept. of
Functional Analysis, Institute of Mathematics, University of {\L}{\'o}d{\'z},
ul. S. Banacha 22, 90-238  {\L}\'od\'z,
 Poland.}
%$^{\cal P}$ Partially supported by {\L}{\'o}d{\'z} University grant 422/93\\
\ftt[5]{Partially
supported by KBN grant 2P 302 21706.}
\ftt[6]{Partially supported by KBN grant 2P 302 08706.}
\setcounter{footnote}{0}
\def\thefootnote{\arabic{footnote}}
\end{titlepage}
\newpage \setcounter{page}{1}

\sec
\section{Introduction}

The inhomogeneous groups $G$ described by means of the semidirect product
$G=H\bowti A$
of a simple Lie group $H$ and the Abelian group $A$ are very important in
physical applications (see e.g. [1]). The most important example is given by
the $D=4$ \poin group $\Pc=SO(3,1)\bowti T^4$ where $H=SO(3,1)$ is the Lorentz
group and $A=T^4$ describes four Abelian translations. The \poin group as a
semidirect product can be supersymmetrized in two ways:
\ben%[
\i[i)] We keep the factor $H=SO(3,1)$ unchanged but we replace the Abelian
subgroup $A=T^4$ by its superextension $T^{4;4}$ , with the $D=4$
superalgebra generators $(\o Q_\ka = (Q_\a))^+$
\bel{1.1}
\ba{rcl}
\{ Q_\a, \o Q_\kb \} &=& 2(\s_\mu)_{\a\kb}P^\mu\,,\\{}
\{ Q_\a, Q_\b \} &=& \{ \o Q_\ka, \o Q_\kb \} = 0\,,\\{}
[Q_\a,P_\mu ] &= & % 0 \,,\qquad
 [\o Q_\ka,P_\mu ] =0 \,,\\{}
[P_\mu,P_\nu] &=& 0\,.
\ea
\ee
The $D=4$ \poin group $\sP$ can be written therefore as (see e.g.\ [2])
\bel{1.2}
\sP=SO(3,1) \bowti T^{4;4}\,,
\ee
where both factors in  \r{1.2} are nonAbelian.
\i[ii)]
One can attach to the factor $H=SO(3,1)$
only two complex supertranslations, e.g.\
generated by $Q_\a$ ($\a=1,2$). In such a way we introduce graded Lorentz
group $\so$ with the Lorentz generators $M_{\mu\nu}$ and two anticommuting
odd
generators $Q_\a$, satisfying the relations
\bel{1.3}
[M_{\mu\nu}, Q_\a] = \frac{i}2 (\s_{\mu\nu})_\a{}^\b Q_\b\,.
\ee
Denoting the graded Abelian group generated by generators ($P_\mu$, $\o
Q_\ka$) by $\o \sT$ one can write
\bl
\bel{1.4a}
\sP = \so \bowti \o \sT
\ee
In fact there is still another possibility, obtained by substituting in $\so$
the generators $Q_\a$ by $\o Q_\ka$ ($\so \rightarrow \o{\so}$) and replacing
in $\o \sT$ the generators $\o Q_\ka$ by $Q_\a$ ($\o \sT \rightarrow \sT$).
One gets
\bel{1.4b}
\sP = \o{\so}\ \bowti \sT\,.
\ee
\el
The advantage of the semidirect products \r{1.4a}-\r{1.4b} is the
preservation of the graded Abelian nature of the second factor. It appears
that this graded Abelian structure will be preserved also after the
$\k$-deformation (compare with [9], where the fourmomenta for \kdef \poin
algebra commute).
\een%]

The main result of this paper is the description of quantum $\k$-deformation
of $D=4$ \poin supergroup in the framework of bicrossproduct Hopf
superalgebras. It appears that for the standard (``bosonic'') quantum groups
the bicrossproduct Hopf algebra $H_1\btc H_2$ provides an attractive
proposal for the general framework of quantum deformations of classical
semidirect product of $H_1$ and $H_2$ supplemented by consistent coalgebra
structure [3 -- 5]. Indeed recently it has been shown that the dual pair of
\kdef \poin algebra [6,7] and \kdef \poin group [8] can be very
well incorporated
[9 -- 11] into the bicrossproduct scheme
\footnote{The general framework proposed recently by Podle\'s and Woronowicz
[12] can be also put into the general bicrossproduct framework.}.
Our aim here is to extend such a framework to the ``super''-case.

Our presentation contains two parts:
\ben%[
\i[i)] It appears that the signature factors which enter into the defining
properties of graded bicrossproduct Hopf algebras are not known in the
literature, and appear in several formulas in a way that is far from obvious.
In Sect.\ 2 we describe the general framework describing graded
bicrossproduct Hopf algebras, extending to super-case the results presented
by Majid [3 -- 5].
\i[ii)] In Sect.\ 3 we describe the $D=4$ \kdef \poin supergroup, presented
firstly in [13], as the graded bicrossproduct Hopf algebra, i.e.\ we provide
the example of our general scheme. In this description the graded Hopf algebra
is the algebra of functions on quantum supergroup with graded set of
generators.
\een%]

It should be stressed that the bicrossproduct description of Hopf algebra $H=
H_1\btc H_2$ implies the bicrossproduct structure of dual Hopf algebra
$\~H = \~H_2 \btc \~H_1$, where the ``tilde'' describes the dual object.
Indeed, the first complete proof of duality between $D=4$ $\k$-\poin quantum
algebra and $D=4$ $\k$-\poin quantum group was obtained in [11] after using
the bicrossproduct form of both dual structures. In the present paper it is
firstly demonstrated that the $D=4$ $\k$-\poin supergroup is an example of
graded bicrossproduct Hopf algebra. Because the bicrossproduct structure of
$D=4$ $\k$-\poin quantum superalgebra is already known [14] we believe that
it is only a technical matter to extend the proof given in [11] to the
supersymmetrized case.

\sec
\section{Bicrossproduct Hopf algebras}
Firstly we shall recall the general definition of bicrossproduct Hopf algebra,
due to Majid [3 -- 5]. Let $H_1$, $H_2$ ($a \in H_1$, $h \in H_2$) are two
Hopf algebras, with the coproducts $\cop(a)=a_\1\tens a_\2$,
$\cop(h)=h_\1\tens
h_\2$. We assume further that
\ben%[
\i[i)] $H_2$ is  a right $H_1$ module algebra, i.e.\ there exists an action
$\a\:H_2\tens H_1 \to{} H_2$ denoted by
\bel{2.1}
\a(a \tens h) = a \trl h
\ee
defining the cross-product $H_1 \bowti H_2$.
\i[ii)] $H_1$ is a left $H_2$-comodule coalgebra, i.e. there exists a
coaction $\b\: H_1 \to{} H_2 \tens H_1$ such that
\bel{2.2}
\b(h)=h^\1 \tens h^\2\,,
\ee
where $h^\1\in H_2$, $h^\2\in H_1$, which defines the crossproduct $H_1 \ptc
H_2$.
\een%]
The following statement is valid:
\begin{tw}\label{t1}(see [3], Theorem 3.3):\\
The linear space $H_1 \tens H_2$ is endowed with the Hopf algebra structure
and defines the right-left bicrossproduct Hopf algebra $H_1 \btc H_2$ if the
action $\a$ and coaction $\b$ satisfy the following compatibility conditions:
\bl
\bel{2.3a}
\e( a \trl h) = \e(a)\e(h)\,,\qquad \b(1) = 1 \tens 1 \,,
\ee
\bel{2.3b}
\cop(a \trl h) = (a_\1 \trl h_\1) h_\2{}^\1 \tens a_\2 \trl h_\2{}^\2\,,
\ee
\bel{2.3c}
\b(hg) = (h^\1 \trl g_\1) g_\2{}^\1 \tens h^\2 g_\2{}^\2\,,
\ee
\bel{2.3d}
h_\1{}^\1(a\trl h_\2) \tens h_\1{}^\2 = (a \trl h_\1) h_\2{}^\1 \tens
h_\2{}^\2\,.
\ee
\el
The multiplication structure in $H_1 \btc H_2$ is defined by
\bel{2.4}
(h \tens a)\cdot (g \tens b) = hg_\1 \tens (a \trl g_\2 ) b
\ee
and the comultiplication looks as follows
\bel{2.5}
\cop(h \tens a)= (h_\1 \tens h_\2{}^\1 a_\1)\tens (h_\2{}^\2 \tens a_\2)\,.
\ee
The antipode and the counit are given by the formulae
\bl
\bel{2.6a}
S(h\tens a) = (1 \tens S(h^\1 a))\cdot(S(h^\2) \tens 1 )\,,
\ee
\bel{2.6b}
\e(h \tens a) = \e(h) \e(a)\,.
\ee
\el
\end{tw}
The simplest examples of bicrossproduct Hopf algebras are:
\ben%[
\i[i)]
Semidirect product of simple Lie algebra and the Abelian algebra, e.g.\
$\Pc=SO(3,1) \btc T^4$, considered  a Hopf algebra with primitive coproduct.
In such a case ($M_{\mu\nu}\in SO(3,1)$, $P_\mu \in T^4$)
\bel{2.7}
\a(P_\mu \tens M_{\rho\t}) = P_\mu \trl M_{\rho\t} = [ P_\mu , M_{\rho\t}]\,,
\qquad
\b(M_{\mu\nu}) = 1 \tens M_{\mu\nu}\,.
\ee
\i[ii)] If $H_1$ is commutative and $H_2$ is cocommutative, the bicrossproduct
$H_1 \btc H_2$ describes the Hopf algebra extension in the sense of Singer
[15].
\i[iii)] The quantum $\k$-deformation of $D=4$ \poin algebra $\U(so(3,1))\btc
\U_\k(T^4)$ [9, 10] as well as the quantum $\k$-deformation od $D=4$ \poin
group $C(SO(3,1)) \ctb C_\k(\~T^4)$ [9 -- 11], where $C(SO(3,1))$ and
$C(\~T^4)$ describe respectively the algebra of functions on the Lorentz
group and the algebra of functions on the Abelian translation group
($C_\k(\~T^4)$ is dual to $\U_\k(T^4)$).
\een%]

\sec
\section{$Z_2$-graded bicrossproduct Hopf superalgebra}
\subsection{$Z_2$-graded Hopf superalgebras}
Let us assume that $\H$ is a $Z_2$-graded algebra, i.e.\ as a vector space
$\H=\H_0 \oplus \H_1$. We define the parity of an element $h\in \H$ as
follows:
\bel{3.1}
p(h)=0 \mbox{ if } h\in\H_0;\qquad p(h)=1 \mbox{ if } h\in\H_1\,.
\ee

We introduce the tensor product of two $Z_2$-graded algebras $\H\tens\H'$ as
$Z_2$-graded algebra with the following multiplication rule
\bel{3.2}
(h\tens h')\cdot(g\tens g') =(-1)^{p(h')p(g)} hg \tens h'g'\,,
\ee
where $h,g\in \H$ and $h',g'\in\H'$. Let us introduce the following definitions
[17 -- 19]:
\begin{df}\label{d1}
The Hopf superalgebra ($\H,\cop,\e,S$) is given by the following axioms:
\ben%[
\i[i)] $\H$ is a \graded algebra,
\i[ii)] The comultiplication map $\cop\: \H \to{} \H\tens \H$ is a
homomorphism of $\H$ and is superassociative, i.e.\
\bel{3.3}
(\cop \tens 1) \cop(h) = (1 \tens \cop) \cop(h)\,,
\ee
where the tensor product is defined by \r{3.2}; besides $\cop(1)=1\tens1$.
\i[iii)] The counit $\e$ is linear map $\H \to{} C$, where
\bel{3.4}
\e(hh')=\e(h)\e(h')\,,\qquad (\e\tens1)\cop(h)=(1\tens\e)\cop(h)=h\,,
\ee
\i[iv)] The antipode $S$ is defined as linear antiisomorphism $\H \to{}\H$
with the property ($m(h\tens h')=h\cdot h'$)
%\bl
\bel{3.5b}
m \circ (S\tens 1)\circ\cop(h) = m \circ (1\tens S)\circ\cop(h) =\e (h)\cdot
1\,.
\ee
\een
\end{df}
If we introduce the graded flip operation:
\bel{3.6}
\s(h\tens h')=(-1)^{p(h)p(h')}h'\tens h\,.
\ee
and the graded opposite coproduct
\bel{3.7}
\cop'(h)=\s \circ \cop(h)
\ee
we obtain also that
\bel{3.8}
(S\tens S)\cop(h) = \cop'(S(h))
\ee
as well as
\bel{3.5a}
S(hh')=(-1)^{p(h)p(h')}S(h')S(h)\,.
\ee
%\een%]
%\end{df}
{\bf Remark 1.} The Hopf superalgebra is a special case of anyonic Hopf
algebras, with $Z_n$-grading, for the special case $n=2$ [20]. The anyonic
algebras are the special case of braided Hopf algebras [21, 22].

\noindent {\bf Remark 2.} The well known examples of noncommutative and
noncocommutative Hopf superalgebras are the quantum deformations  of universal
enveloping algebras of simple Lie superalgebras (see e.g.\ [17 -- 19]) as
well as the quantum deformation of the algebra of functions on a simple Lie
supergroups (see e.g.\ [23]).

\subsection{Bicrossproduct Hopf superalgebras}
Let us assume that $\H_1$, $\H_2$ are two Hopf superalgebras (we put
$h,g,f\in \H_1$, $a,b,c\in\H_2$). We assume further that
\ben%[
\i[i)] $\H_2$ is a right $\H_1$ module, with the action $\a$ (see \r{2.1}),
satisfying the grading property
\bl
\bel{3.9a}
ab \trl h = (-1)^{p(h_\1)p(b)} (a \trl h_\1)(b \trl h_\2)
\ee
\bel{3.9b}
a \trl (hg) = (a \trl h) \trl g
\ee
\el
\i[ii)] $\H_1$ is a left $\H_2$-comodule cosuperalgebra with the action $\b$
(see \r{2.2}), satisfying $\b(1)=1\tens 1$ and the following properties
\bl
\bel{3.10a}
(1 \tens \b) \circ \b = (\Delta \tens 1) \circ \b\,,
\ee
\bel{3.10b}
(\e \tens 1) \circ \b(h) = 1_{\H_1} \tens h\,,
\ee
\bel{3.10c}
(1 \tens \cop) \b(h) = m_{12}\s_{23} (\b\tens\b)\cop
\ee
\el
\een%]
\begin{tw}\label{t2}
The linear \graded space (superspace) $\H_1\tens\H_2$ is endowed with the
Hopf superalgebra structure and defines the \graded bicrossproduct Hopf
superalgebra $\H_1\btc\H_2$ with the following definitions of multiplication,
comultiplication, counit and antipode
\bel{3.10}
(h\tens a)(g\tens b) =(-1)^{p(a)p(g_\1)}h g_\1 \tens (a \trl g_\2) b\,,
\ee
\bel{3.11}
\cop(h\tens a) = (-1)^{p(h_\2{}^\2)p(a_\1)}
           h_\1\tens h_\2{}^\1 a_\1 \tens h_\2{}^\2 \tens a_\1\,,
\ee
\bel{3.12}
\e(h\tens a) =\e(h) \cdot \e (a)\,,
\ee
\bel{3.13}
S(h\tens a) = (-1) ^{p(h^\2)[p(h^\1)+p(a)]}(1 \tens S(h^\1 a))(S(h^\2)\tens 1)
\ee
if the following compatibility conditions are satisfied (compare with
\r{2.3a}-\r{2.3d}):
\bl
\bel{3.14a}
\e(a \trl h) = \e(a) \e(h)\,,
\ee
\bel{3.14b}
\cop(a \trl h) (-1)^{p(a_\2)[p(h_\1)+p(h_\2{}^\1)]}
(a_\1\trl h_\1 ) h_\2{}^\1 \tens a_\2 \trl h_\2{}^\2 \,,
\ee
\bel{3.14c}
\b(hg)=(-1)^{p(h^\2)[p(g_\1)+p(g_\2{}^\1)]}
(h^\1 \trl g_\1)g_\2{}^\1 \tens h^\2 g_\2{}^\2\,,
\ee
\bel{3.14d}
h_\1{}^\1 (a \trl h_\2) \tens h_\1{}^\2 =
(-1)^{p(a)p(h_\2{}^\1)+p(h_\1)p(h_\2{}^\2)}
(a \trl  h_\1) h_\2{}^\1 \tens h_\2{}^\2\,.
\ee
\el
\end{tw}

In the proof we check directly that with the definitions \r{3.10}-\r{3.13}
 and the properties \r{3.14a}-\r{3.14d} the axioms of the graded
bicrossproduct $\H_1 \btc \H_2$ being a Hopf superalgebra are satisfied. In
particular we have to check the associativity of the multiplication \r{3.10},
i.e.\
\bel{3.15}
[(h\tens a)(g \tens b)] (f \tens c) = (h\tens a)[(g \tens b) (f \tens c) ]\,,
\ee
the coassociativity of the coproduct \r{3.11} (see \r{3.3}) and, what is the
most complicated part, the homomorphism property of the coproduct
\bel{3.16}
\cop((h\tens a)(g \tens b)) = \cop (h \tens a) \cop (g \tens b)\,.
\ee
The graded bicrossproducts described in this sections can be used for the
description of quantum deformations of inhomogeneous supergroups.

\sec
\section{An example: $D=4$ \kdef \poin supergroup in bicrossproduct form}

The $D=4$ \kdef \poin supergroup has been obtained in [13] by quantization of
the $r$-matrix Poisson bracket, with the following choice of the classical
$r$-matrix for $D=4$ \poin superalgebra (see also [24])
\bel{4.1}
r = N_i \wedge P_i - \frac{i}{4} Q_\a \wedge \o Q_\ka\,,
\ee
where $M_{\mu\nu}= (M_i,N_i)$ are Lorentz generators, $P_\mu$ describe the
fourmomenta and $Q_{\alpha}$, $\o Q_\ka$ are four supercharges described by the
doublet of Weyl $2$-spinors. The Hopf superalgebra of the functions on $D=4$
\kdef \poin supergroup is described by the following relations (see [13]):
\ben%[
\i[a)] Algebra
% Wstawka z KLMS
%\ben%[
\item[i)] Lorentz sector $(A_{\alpha}^{\beta},
A_{\dot\alpha}^{\dot\beta})$ (we use the spinorial representation of the
Lorentz generators, e.g.\
$ L_i = {1\over 4}(\sigma_i)_{\alpha}^{\beta}L^{\alpha}_{\beta}
+(\overline\sigma_i)_{\dot\alpha}^{\dot\beta} L^{\dot\alpha}_{\dot\beta} $):

The Lorentz subgroup parameters are classical, i.e.
\bel{4.2}
[A_{\alpha}^{\beta}, A_{\gamma}^{\delta}] = [A_{\alpha}^{\beta},
A_{\dot\gamma}^{\dot\delta} ] = [A_{\dot\alpha}^{\dot\beta},
A_{\dot\gamma}^{\dot\beta} ] = 0
\ee
\item[ii)] Translations $(X_{\mu})$ (we denote $\theta =
{{\theta_1}\choose{\theta_2}}, \overline\theta = {{\theta_{\dot
1}}\choose{\theta_{\dot 2}}}$):
\bel{4.3}
\begin{array}{ll}
[X^i, X^j] = {i\over {8\kappa}}\theta^T
\sigma^i (\hbox{\bf 1}_2 - (AA^+)^{-1})\sigma^j
\overline\theta-{i\over {8\kappa}}\theta^T \sigma^j (\hbox{\bf
1}-(AA^+)^{-1})\sigma^i \overline\theta&\null\\ \\{}
[X^0, X^j] = -{i\over {\kappa}}X^j + {i\over
{8\kappa}}\theta^T \lbrack \sigma^j,
(AA^+)^{-1}\rbrack \overline\theta&\null
\end{array}
\ee
\bel{}
\begin{array}{lll}
[ A_{\alpha}^{\beta}, X^i]&=&{1\over
{2\kappa}}((A\sigma_{n})_{\alpha}^{\beta} \Lambda^i_n (A) -
(\sigma^i\cdot\!A)_{\alpha}^{\beta})\\ \\
{}[ A_{\alpha}^{\beta}, X^0]&=&{1\over {2\kappa}}(A\sigma_i
)_{\alpha}^{\beta} \Lambda^0_i (A)
\end{array}
\ee
\item[iii)] Supertranslations
\bel{4.5}
\{ \theta^{\alpha}, \theta^{\beta}\}=\{ \theta^{\dot\alpha},
\theta^{\dot\beta}\} = 0\quad
\{ \theta^{\alpha}, \theta^{\dot\beta}\}={i\over {2\kappa}}(\hbox{\bf 1} -
(AA^+)^{-1})^{\dot\beta \alpha}
\ee
\bel{4.6}
\begin{array}{lll}
\{ X^i, \theta_{\alpha} \} &=&{1\over
{4\kappa}}(\theta^T\sigma^i)_{\gamma} (\hbox{\bf 1}_2 -
(AA^+)^{-1})^{\gamma}_{\alpha}\\ \\
\{ X^0, \theta_{\alpha}\}&=&-{1\over
{4\kappa}}\theta^T_{\gamma}(\hbox{\bf 1}_2 +
(AA^+)^{-1})^{\gamma}_{\alpha}
\end{array}
\ee

\bel{4.7}
\begin{array}{lll}
\{ A_{\alpha}^{\beta}, \theta^{\gamma}\}&=&\{
A_{\dot\alpha}^{\dot\beta}, \theta^{\gamma} \} = 0
\end{array}
\ee
%\een
\i[b)] Coalgebra \\
\bel{4.8}
\begin{array}{ll}
\Delta (X_{\mu})&=X_{\mu}\otimes\hbox{\bf 1} +
\Lambda_{\mu}^{\nu}(A)\otimes X_{\nu}-
{i\over 2}({A^{-1}_{\alpha}}^{\beta} \sigma^{\mu}_{\beta\dot\gamma}
\theta^{\dot\gamma} \otimes \theta^{\alpha} +
\theta^{\alpha}\sigma^{\mu}_{\alpha\dot\beta}A^{-1\beta}_{\dot\gamma}
\otimes \theta^{\dot\gamma})\\
\Delta (\theta_{\alpha})&=\theta_{\alpha}\otimes\hbox{\bf
1}+(A^{-1})^{\beta}_{\alpha}\otimes\theta_{\beta}\\
\Delta (A_{\alpha}^{\beta})&=A_{\alpha}^{\gamma}\otimes
A_{\gamma}^{\beta}
\end{array}
\ee
\i[c)] Antipodes:
\bel{4.9}
S(X^{\mu})=-\Lambda^{\mu}_{\nu}(A^{-1})X^{\nu}\quad
S(A_{\alpha}^{\beta})=(A^{-1})_{\alpha}^{\beta}
\ee
$$
S(\theta^{\alpha})=-A_{\beta}^{\gamma}\theta^{\beta}
$$
\een
In such a way we have obtained the complete set of relations describing
the $\kappa$-deformation of $N=1$ Poincar{\'e} supergroup.
%koniec wstawki

In order to put the Hopf superalgebra \r{4.2}-\r{4.9} in the bicrossproduct
algebra form we should introduce the complexified chiral superspace
coordinates (see e.g. [25]) as follows
\bel{4.10}
z_\mu = x_\mu + \frac{i}2 \theta _\a (\s_\mu)^{\a\kb} \o \theta _\kb\,.
\ee

Let us introduce the following two \kdef Hopf superalgebras:
\ben%[
\i[i)] The algebra of functions $C(z_\mu,\theta_\a)$ on \kdef chiral
superspace $(z_\mu,\theta_\a)$, with the following Hopf superalgebra
relations:
\bel{4.11}
\ba{rclrcl}
[z_i,z_j]&=&0\,,\qquad & [z_0,z_i]&=& -\frac i \k z_i\,,\\
{}[z_0,\theta_\a] &=& -\frac i {2\k} \theta_\a\,,\quad& [z_i,\theta_\a]&=&0,
\qquad \{\theta_{\alpha}, \theta_{\beta}\}=0.
\ea
\ee
and
\bel{4.12}
\ba{rcl}
\cop(z_\mu) &=&z_\mu\tens 1 + 1 \tens z_\mu\,,\\
\cop(\theta_\a)&=&\theta_\a \tens 1 + 1 \tens \theta_\a\,.
\ea
\ee
\i[ii)] The classical Hopf superalgebra of functions $C(A_{\a\b}, A_{\kg\kb},
\o \theta_\ka)$ on the superextension of classical Lorentz group, with the
following defining relations:
\bel{4.13}
\ba{rclllll}
[A_{\a\b},A_{\g\d}]&=& [A_{\a\b}, A_{\kg\kd}]&=& [A_{\ka\kb},A_{\kg\kd}] &=&
0\,,\\
{}[A_{\a\b}, \o \theta_\kg]&=& [A_{\ka\kb}, \o \theta_\kg] &=& 0\,,\\
\{\o\theta_\ka, \o \theta_\kb\} &=& 0
\ea
\ee
and
\bel{4.14}
\ba{l}
\cop A_{\a\b} = A _{\a\g} \tens A_{\g\b}\,,\qquad \cop A_{\ka\kb} =
A_{\ka\kg}\tens A_{\kg\kb}\,,\\
\cop\o \theta _\ka = \o \theta _\ka \tens 1 + (A_{\ka\kb})^{-1} \tens \o
\theta_\kb\,.
\ea
\ee
\een%]
One can prove the following statement:
\begin{tw}\label{t3}
The $D=4$ \kdef \poin supergroup can be described as the graded
bicrossproduct Hopf superalgebra
\bel{4.15}
C_\k(\sP) = C(z_\mu, \theta_\a) \btc C(A_{\a\b}, A_{\ka\kb}, \o \theta_\ka)\,,
\ee
with the following definition of the action $\a$
\bel{4.16}
\ba{rcl}
\ota \trl z_i &=& -\frac i {2\k} \[ 1- (A^+A)^{-1} \]_\ka{}^\kb \otb\,,\\
\ota \trl z_0 &=& -\frac i {2\k} (A^+A)_\ka{}^\kb \otb\,,\\
\ota \trl \theta_\b &=& -\frac i {2\k} \[ 1- (AA^+)^{-1} \]_{\ka\b} \,,\\
A_{\a\b} \trl z_i &=& \frac1{2\k}(A\s_k)_{\a\b} \Lambda_{ik} (A,\o A)-(\s_i
A)_{\a\b})\,,\\
\o A_{\ka\kb} \trl z_i &=& \frac1{2\k}[(\s_i \o A)_{\ka\kb} \Lambda_{lk}
(A,\o A)-(\o A \s_i)_{\ka\kb}]\,,\\
A_{\a\b} \trl z_0 &=& \frac1{2\k}(A\s_i)_\a{}^\b \Lambda_{i0} (A,\o  A)\,,\\
\o A_{\ka\kb} \trl z_0 &=& \frac1{2\k}(\s_i\o A)_\ka{}^\kb \Lambda_{i0} (A,\o
 A)\,,\\
A_{\a\b} \trl \theta_\g &=& A_{\ka\kb} \trl \theta_\g =0 \\
\ea
\ee
and coaction $\b$
\bel{4.17}
\ba{rcl}
\b(z_\mu)&=& \Lambda_\mu{}^\nu (A,\o A) \tens z_\nu - i
(A^{-1)}\s_\mu)_{\a\kb} \otb \tens \theta_\a\,,\\
\b(\theta_\a) &=& (A^{-1})_{\b\a} \tens \theta_\b\,.
\ea
\ee
\end{tw}
The proof is obtained by performing the transformation \r{4.10} of the
superalgebra basis \r{4.2}-\r{4.9} and showing that the resulting Hopf
superalgebra fits into the graded bicrossproduct superalgebra framework.

\sec
\section{Final Remarks}
In this paper we give not known in the literature general definition of
graded bicrossproduct Hopf superalgebra and we provided as well an example.
Such a scheme can be extended to the case of anyonic groups [20] and braided
groups [21, 22]. To our knowledge only the crossproducts of braided
cocommutative Hopf algebra with the quasitriangular Hopf algebras have been
considered [26]. It is an interesting task to introduce the general
braided bicrossproduct of braided Hopf algebras and provide some
nontrivial examples.

\newpage
\def\bib{\bibitem}

\end{document}